\documentclass[graybox, envcountchap]{svmult}

\usepackage{amsmath,amsfonts, amssymb}

\usepackage{mathptmx}    
\usepackage{helvet}           
\usepackage{courier}         

\usepackage{makeidx}      
\usepackage{graphicx}      

\usepackage{multicol}      
\usepackage[bottom]{footmisc} 

\makeindex            

\newcommand{\field}[1]{\mathbb{#1}}

\newcommand{\R}{\field{R}}

\begin{document}

\title{High-order ADI schemes for convection-diffusion equations with mixed derivative terms}
\titlerunning{High-order ADI schemes} 
\author{B. D{\"u}ring, M. Fourni{\'e} and A. Rigal}
\institute{B.  D{\"u}ring \at  Department of Mathematics, University of Sussex, Pevensey 2, 
 Brighton, BN1 9QH, United Kingdom, \email{b.during@sussex.ac.uk} 
\and  M. Fourni{\'e} \at Institut de Math{\'e}matiques de Toulouse,
 Equipe `Math{\'e}matiques pour l'Industrie et la Physique',
  CNRS, Unit\'e Mixte 5219,  Universit{\'e}s de Toulouse,
 118, route de Narbonne, 31062 Toulouse Cedex, France \email{michel.fournie@math.univ-toulouse.fr} \and A. Rigal \at Institut de Math{\'e}matiques de Toulouse,
 Equipe `Math{\'e}matiques pour l'Industrie et la Physique',
  CNRS, Unit\'e Mixte 5219,  Universit{\'e}s de Toulouse,
 118, route de Narbonne, 31062 Toulouse Cedex, France \email{alain.rigal@math.univ-toulouse.fr}}

\maketitle

\abstract*{We present new high-order Alternating Direction Implicit (ADI) schemes
for the numerical 
solution of initial-boundary value problems for convection-diffusion equations with mixed derivative terms.
Our approach is based on the unconditionally stable ADI scheme 
proposed by Hundsdorfer \cite{Hund02}. 
Different numerical discretizations which lead
to schemes which are fourth-order accurate in space and second-order accurate in time are discussed.
}

\abstract{We present new high-order Alternating Direction Implicit (ADI) schemes
for the numerical 
solution of initial-boundary value problems for convection-diffusion equations with cross derivative terms.
Our approach is based on the unconditionally stable ADI scheme 
proposed by Hundsdorfer \cite{Hund02}. 
Different numerical discretizations which lead
to schemes which are fourth-order accurate in space and second-order accurate in time are discussed.
}

\section{Introduction}
\label{sec:1}
We consider the multi-dimensional convection-diffusion equation
\begin{equation}
\label{convDiff}
u_t =\mathrm{div} (D\nabla u) + c \cdot \nabla u 
\end{equation}
on a rectangular domain $\Omega \subset  \mathbb{R}^2$, supplemented
with initial and boundary conditions. In (\ref{convDiff}), 
$$c=\left ( \begin{array}{c} 
c_1\\
c_2\\
\end{array}
 \right ),\quad D=\left ( 
\begin{array}{cc} 
d_{11}&d_{12}\\
d_{21}&d_{22}\\
\end{array}
\right ),
$$ 
are a given nonzero 
convection vector and a given, fully populated (non-diagonal), and
positive definite diffusion matrix, respectively.
Thus, both mixed derivative and convection terms are present in (\ref{convDiff}).

 After rearranging, problem \eqref{convDiff} may be formulated as
\begin{equation}
\label{convDiff2}
\frac{\partial u(x,y,t)}{\partial t} =  \underbrace{ (d_{12}+d_{21}) \frac{\partial^2 u}{\partial x \partial y}}_{=:F_0(u)} + 
\underbrace{(c_1 \frac{\partial u}{\partial x} + d_{11}\frac{\partial^2 u}{\partial x^2})}_{=:F_1(u)} 
+ \underbrace{(c_2 \frac{\partial u}{\partial y} + d_{22}\frac{\partial^2 u}{\partial y^2})}_{=:F_2(u)}.
\end{equation}

This type of convection-diffusion equations with mixed derivatives
arise frequently in many applications, e.g.\ in financial mathematics
for option pricing in stochastic volatility models or in numerical mathematics
when coordinate transformations are applied.
Such transformations are particularly useful to allow working on
simple (rectangular) domains or on uniform grids (to have better accuracy).
Thus, this approach allows to consider complex domains or 
to define non-uniform meshes to take into account the stiffness 
behavior of the solution in some part of the domain. 

In the mathematical literature, there exist a number of
numerical approaches to approximate solutions to \eqref{convDiff},
e.g.\ finite difference schemes, spectral methods, finite volume and
finite element methods. Here, we consider \eqref{convDiff} on a
rectangular domain $\Omega \subset \mathbb{R}^2$. In this situation a
finite difference approach seems most straight-forward. 

The Alternating Direction Implicit (ADI) method introduced by
Peaceman and Rachford \cite{PeacmanRachford95}, Douglas \cite{Douglas,
Gunn}, Fairweather and Mitchell \cite{Fairweather}
is a very powerful method that is especially
useful for solving parabolic equations on rectangular domains. Beam
and Warming \cite{BeamWarming}, however, have shown that no simple 
ADI scheme involving only discrete solutions at time levels $n$ and $n+1$ can be
second-order accurate in time in the presence of mixed derivatives
($F_0\neq 0$ in \eqref{convDiff2}). To overcome this limitation and construct an
unconditionally stable ADI scheme of second order in time, a number of
results have been given by Hundsdorfer \cite{Hund02, Verwer} and 
more recently by in 't Hout and Welfert \cite{HoutWelfert07}.
These schemes are second-order accurate in time and space.

High-Order Compact (HOC) schemes (see, e.g.\ \cite{Gupta,Rigal99}) employ a
nine-point computational stencil using the eight 
neighbouring points of the reference grid point only and show good
numerical properties.
 Several papers consider the application of HOC
schemes (fourth order accurate in space) for two-dimensional convection-diffusion problems {\em with mixed
derivatives} \cite{FournieKaraa06,DuringFournie12} but {\em without ADI} splitting. 
Moreover, the HOC approach introduces a high algebraic complexity in the derivation of
the scheme.

We are interested in obtaining efficient, {\em high-order\/}
ADI schemes, i.e.\ schemes which have a consistency order equal to two in time
and to four in space, which are unconditionally stable and robust (no oscillations).
We combine the second-order ADI splitting scheme presented in 
\cite{Hund02,HoutWelfert07} with different high-order schemes to
approximate $F_0, F_1, F_2$ in \eqref{convDiff2}.
We note that some results on coupling {\em HOC with ADI\/} have been
presented in \cite{KaraaZhang04}, however, {\em without 
mixed derivative terms} present in the equation. 

Up to the knowledge of the authors there are currently no results
for ADI-HOC in the presence of mixed derivative terms.
In this preparatory work we validate the coupling of ADI and HOC by numerical experiments.


\section{Splitting in time}

In time, we consider the following splitting scheme presented in 
\cite{Hund02,HoutWelfert07}. We consider \eqref{convDiff2}, and 
we look for a (semi-discrete) approximation $U^n \approx u(t_n)$
with $t_n =n\Delta_t$ for a time
step $\Delta_t$.
The scheme used corresponds to 
\begin{equation}
\label{eq:HVscheme}
\left \{
\begin{array}{l}
Y^0 = U^{n-1} +\Delta_t F(U^{n-1}),\\
Y^1 = Y^{0} +\theta \Delta_t(F_1 (Y^1 )-F_1 (U^{n-1})),\\
Y^2 = Y^{1} +\theta \Delta_t(F_2 (Y^2 )-F_2 (U^{n-1})),\\
\tilde{Y}^0 = Y^0 + \sigma \Delta_t(F(Y^2)-F(U^{n-1})),\\
\tilde{Y}^1 = \tilde{Y}^{0} + \theta \Delta_t(F_1 (\tilde{Y}^1 )-F_1 (Y^{2})),\\
\tilde{Y}^2 = \tilde{Y}^{1} + \theta \Delta_t(F_2 (\tilde{Y}^2 )-F_2 (Y^{2})),\\
U^n = \tilde{Y}^2,
\end{array}
\right .
\end{equation}
with constant parameters $\theta$ and $\sigma,$ and $F = F_0+ F_1+F_2.$
To ensure second-order consistency in time we choose $\sigma = 1/2$.
The parameter $\theta$ is arbitrary and typically fixed to $\theta=1/2$. The choice
of $\theta$ is discussed in \cite{Hund02}. 
Larger $\theta$ gives stronger damping of implicit terms and  lower
values return better accuracy (some numerical results for $\theta = 1/2+\sqrt{3}/6$ are given in section~\ref{sec:numexp}).

We note that $F_0$ is treated explicitly,  whereas $F_1,F_2$
(unidirectional contributions in $F$)  are treated implicitly. 
In the following section, we discuss different high-order (fourth order) strategies for the
discretization in space.


\section{High-order approximation in space}

For the discretization in space, we replace the rectangular domain $\Omega=[L_1,R_1]
\times [L_2,R_2]\subset \R^2$ with $R_1>L_1$, $R_2>L_2$ by a uniform grid $Z=\{
x_{i}\in \left[L_1,R_1%
 \right]:$ $x_{i}=L_1+(i-1)\Delta_x$, $i=1,\dots,N\}\times\{ y_{j}\in \left[L_2,R_2%
 \right]:$ $y_{j}=L_2+(j-1)\Delta_y$, $j=1,\dots,M\}$ consisting of
 $N \times M$ 
 grid points, with space steps $\Delta_x=(R_1-L_1)/(N-1)$ and
 $\Delta_y=(R_2-L_2)/(M-1)$. Let $u_{i,j}$ 
 denote the approximate solution in $(x_{i},y_j)$ at some fixed time (we omit the superscript $n$ to simplify the 
 notation).

We present different fourth-order schemes to approximate $F_0, F_1, F_2$ in \eqref{eq:HVscheme}.
The first one uses five nodes in each direction and the second one is
compact.
Both schemes are considered with boundary conditions of either
periodic or Dirichlet type.

\subsection{Fourth-order scheme using five nodes}
\label{sec:HO5}

We denote by $\delta_{x0}$,  $\delta_{x+}$ and
$\delta_{x-}$,  the standard central, forward and backward finite
difference operators, respectively. The second-order central
difference operator is denoted by $\delta_x^2$,
$$
\delta_x^2 u_{i,j}=\delta_{x+}\delta_{x-}
u_{i,j}=\frac{u_{i+1,j}-2u_{i,j}+u_{i-1,j}}{\Delta_x^2}.
$$ 
The difference operators in the $y$-direction, $\delta_{y0}$,  $\delta_{y+}$, $\delta_{y-}$
and $\delta_y^2$, are defined analogously.
Then it is possible to define fourth-order approximations based on,
\begin{equation}
\hspace*{-0.7cm}
\label{eq:4thFDS}
\left .
\begin{array}{ll}
\displaystyle  (u_{x})_{i,j} & \displaystyle \approx \biggl ( 1 - \frac{\Delta_x^2}{6}\delta_x^2 \biggr ) \delta_{x0} u_{i,j}
= \frac{-u_{i+2,j}+8u_{i+1,j}-8u_{i-1,j}+u_{i-2,j}}{12\Delta_x} , \\
\displaystyle (u_{y})_{i,j}&\displaystyle  \approx \biggl ( 1 - \frac{\Delta_y^2}{6}\delta_y^2 \biggr ) \delta_{y0} u_{i,j}
= \frac{-u_{i,j+2}+8u_{i,j+1}-8u_{i,j-1}+u_{i,j-2}}{12\Delta_y},\\
\displaystyle (u_{xx})_{i,j}& \displaystyle \approx 
 \biggl ( 1 - \frac{\Delta_x^2}{12}\delta_x^2 \biggr ) \delta_{x}^2 u_{i,j} = \frac {  - u_{i+2,j}+16u_{i+1,j}-30u_{i,j}+16u_{i-1,j} -u_{i-2,j}}{12\Delta_x^2},\\
 \displaystyle (u_{yy})_{i,j} & \displaystyle \approx
 \biggl ( 1 - \frac{\Delta_y^2}{12}\delta_y^2 \biggr ) \delta_{y}^2 u_{i,j} 
 = \frac{ -u_{i,j+2}+ 16u_{i,j+1}-30u_{i,j}+16u_{i,j-1} -u_{i,j-2} }{12\Delta_y^2},\\
\displaystyle  (u_{xy})_{i,j} &\displaystyle \approx \biggl ( 1 - \frac{\Delta_x^2}{6}\delta_x^2
 \biggr ) \delta_{x0} \biggl ( 1 - \frac{\Delta_y^2}{6}\delta_y^2
 \biggr ) \delta_{y0} u_{i,j}\\
&\displaystyle = \frac{1}{144\Delta_x\Delta_y} \bigl [ 64(
  u_{i+1,j+1}-u_{i-1,j+1}+u_{i-1,j-1} -u_{i+1,j-1})\\
&\displaystyle \hspace*{2cm}+8(-u_{i+2,j+1}-u_{i+1,j+2}+u_{i-1,j+1}+u_{i-2,j+1}\\
&\displaystyle \hspace*{2.6cm}-u_{i-2,j-1}-u_{i-1,j-2}+u_{i+1,j-2} +u_{i+2,j-1})\\
&\displaystyle \hspace*{2cm} +( u_{i+2,j+2}-u_{i-2,j+2}+u_{i-2,j-2} -u_{i+2,j-2})\bigr ].
\end{array}
\right .
\end{equation}
For each differential operators appearing in  $F_0$, $F_1$ and $F_2$,
we use these five-points fourth-order difference formulae.

Combining this spatial discretization with the time splitting
\eqref{eq:HVscheme}, we obtain a high-order, five-points ADI scheme denoted
HO5. Its order of consistency is two in time and four in space.


\subsection{Fourth-order compact scheme}
\label{sec:HOC}

We start by deriving a fourth-order HOC scheme for
\begin{equation}
\label{eq:F1}
F_1(u)=d_{11}\frac{\partial^2u}{\partial x^2}+c_{1} \frac{\partial
  u}{\partial x} = g, 
\end{equation}
with some arbitrary right hand side $g.$ 
We employ central difference operators to approximate the
derivatives in \eqref{eq:F1} using
\begin{align}
\label{eq:ux}
\frac{\partial u}{\partial x}(x_i,y_j) &= \delta_{x0} u_{i,j}
-\frac{\Delta_x^2}{6} { \frac{\partial^3 u}{\partial x^3}(x_i,y_j)} +
\mathcal{O}(\Delta_x^4),\\
\label{eq:uxx}
\frac{\partial^2 u}{\partial x^2}(x_i,y_j) &= \delta_{x}^2 u_{i,j}
-\frac{\Delta_x^2}{12} { \frac{\partial^4 u}{\partial x^4}(x_i,y_j)} +
\mathcal{O}(\Delta_x^4).
\end{align}
By differentiating \eqref{eq:F1}, we can compute the following auxiliary
relations for the derivatives appearing in \eqref{eq:ux}, \eqref{eq:uxx} (in the following, for the sake of brevity we omit the argument
$(x_i,y_j)$ of the continuous functions)
\begin{align}
\label{eq:uxxx}
  \frac{\partial^3u}{\partial x^3}&=\frac{1}{d_{11}} \frac{\partial g}{\partial x} 
-\frac{c_1}{d_{11}} \frac{\partial^2u}{\partial x^2},\\
\label{eq:uxxxx}
\frac{\partial^4u}{\partial x^4}&=\frac{1}{d_{11}} \frac{\partial^2 g}{\partial x^2} 
-\frac{c_1}{d_{11}} { \frac{\partial^3u}{\partial x^3}}=\frac{1}{d_{11}} \frac{\partial^2 g}{\partial x^2} 
-\frac{c_1}{d_{11}} \left ( \frac{1}{d_{11}}\frac{\partial g}{\partial x} - \frac{c_1}{d_{11}}\frac{\partial ^2 u}{\partial x^2} \right ).
\end{align}
Hence, using \eqref{eq:uxxx} and 
\eqref{eq:uxxxx} in \eqref{eq:ux} and \eqref{eq:uxx}, respectively,
equation \eqref{eq:F1} can be approximated by
\begin{equation}
\label{eq:HOCpre}
  d_{11}\delta_{x}^2 u_{i,j} + c_{1} \delta_{x0} u_{i,j}= g_{i,j} +\frac{\Delta_x^2}{12} \left ( \frac{c_1}{d_{11}} \frac{\partial g}{\partial x} + \frac{\partial^2 g}{\partial x^2}  -
        \frac{c_1^2}{d_{11}} \frac{\partial^2 u}{\partial x^2} \right ) +
\mathcal{O}(\Delta_x^4).
\end{equation}
We note that all derivatives on the right hand side of \eqref{eq:HOCpre} can be
approximated on a compact stencil using second-order central difference operators. 
This yields a high-order compact scheme of fourth order for
\eqref{eq:F1} which is given by
\begin{equation}
\label{eq:HOC-F1}
d_{11}\delta_{x}^2 u_{i,j}+ c_{1} \delta_{x0} u_{i,j}   +  \frac{\Delta_x^2}{12}   \frac{c_1^2}{d_{11}} \delta_{x}^2  u_{i,j}  = g_{i,j} +  \frac{\Delta_x^2}{12} 
\left ( \frac{c_1}{d_{11}} \delta_{x0} g_{i,j}  + \delta_{x}^2  g_{i,j}
\right ).
\end{equation}
In a similar fashion we can discretize the operator $F_2(u)=g$ by
a high-order compact scheme of fourth order given by
\begin{equation}
\label{eq:HOC-F2}
d_{22}\delta_{y}^2 u_{i,j}+ c_{2} \delta_{y0} u_{i,j}   +
\frac{\Delta_y^2}{12}   \frac{c_2^2}{d_{22}} \delta_{y}^2  u_{i,j}  = 
g_{i,j} +  \frac{\Delta_y^2}{12} 
\left ( \frac{c_2}{d_{22}} \delta_{y0} g_{i,j}  + \delta_{y}^2  g_{i,j}
\right ).
\end{equation}
Defining vectors $U=(u_{1,1},\dots,u_{N,M})$ and
$G=(g_{1,1},\dots,g_{N,M})$, we can state these schemes \eqref{eq:HOC-F1} and \eqref{eq:HOC-F2} in matrix form
$A_x U= B_x G$ (for $F_1(u)=g$) and $A_y U= B_y G$ (for $F_2(u)=g$), respectively.
We apply these HOC schemes to find the unidirectional contributions $Y^1$,
$\tilde{Y^1}$, and $Y^2$, $\tilde{Y^2}$ in
\eqref{eq:HVscheme}, respectively.
For example, to compute 
 $$Y^1 = Y^{0} +\frac {\Delta_t}2 (F_1 (Y^1 )-F_1 (U^{n-1}))$$
in the second step of \eqref{eq:HVscheme} (which is equivalent to  
$F_1(Y^1-U^{n-1})=-\frac{2}{\Delta_t}(Y^0-Y^1)$), 
we use $A_x(Y^1-U^{n-1})=B_x(-\frac{2}{\Delta_t}(Y^0-Y^1))$ that can be rewrite into 
$$\left (B_x - \frac {\Delta_t}2  A_x \right )Y^1 = B_x Y^0 -\frac {\Delta_t}2  A_x U^{n-1}.$$
Note that the matrix $(B_x - (\Delta_t/2) \,A_x)$ appears twice in
\eqref{eq:HVscheme}, in steps
two and five. Similarly, $(B_y -(\Delta_t/2)\, A_y)$ appears in
steps three and six of \eqref{eq:HVscheme}. Hence, using LU-factorisation,
only two matrix inversions are necessary in each time step of \eqref{eq:HVscheme}.
Moreover, for the case of constant coefficients, these matrices can be
LU-factorized before iterating in time to obtain an even more efficient algorithm.  

To compute $Y^0$ and $\tilde{Y^0}$ in steps one and four of
\eqref{eq:HVscheme} which require evaluation of $F_0$ (mixed term) we use
an explicit approximation using the five-points fourth-order formulae \eqref{eq:4thFDS}.

Combining this spatial discretization with the time splitting
\eqref{eq:HVscheme}, we obtain a high-order compact ADI scheme denoted HOC. Its order of
consistency is two in time and four in space.  

\section{Numerical experiments}
\label{sec:numexp}

We present numerical experiments on a square domain $\Omega =[0,1]\times [0,1]$ for two types of boundary conditions, periodic and Dirichlet type.
The initial condition is given at time $T_0 = 0$ and the solution is computed at the final time $T_f = 0.1$ with different meshes  
$\Delta_x=\Delta_y=h$ and different time steps $\Delta_t$.
In our numerical tests we focus on the errors with respect to time and to space.  

In the first part, we consider the periodic boundary value problem
considered in \cite{HoutWelfert07}.
We implement the scheme detailed in \cite{HoutWelfert07} based on
second-order finite difference approximations (referred to as CDS below)
and compare its behaviour to our new schemes HO5
(section~\ref{sec:HO5}) and HOC (section~\ref{sec:HOC}). 
In the second part, we consider Dirichlet boundary conditions and
restrict our study to the more interesting HOC scheme.
In that part, we extend the splitting scheme to a convection-diffusion
equation with source term. 

\subsection{Periodic boundary conditions}

The problem given in \cite{HoutWelfert07} is formulated on the domain  $\Omega = [0,1]\times[0,1]$. The solution $u$ satisfies
\eqref{convDiff}
 where
$$
c=-\left ( \begin{array}{c} 
2\\
3\\
\end{array}
 \right ),\quad D=0.025 \left ( 
\begin{array}{cc} 
1&2\\
2&4\\
\end{array}
\right ),
$$
with periodic boundary conditions and initial condition $u(x,y,T_0) =e^{-4(\sin^2(\pi x) + \cos^2(\pi y))}$.
We employ the splitting (\ref{eq:HVscheme}) with $\sigma=1/2$ and $\theta=1/2$.

We first present a numerical study to compute the order of convergence
in time of the schemes CDS, HO5 and HOC.
Asymptotically, we expect the error $\epsilon$ to converge as 
$$\epsilon = C \Delta_t^m$$
at some rate $m$ with $C$ representing a constant. This implies 
$$\log(\epsilon) = \log(C) + m \log(\Delta_t).$$
Hence, the double-logarithmic plot $\epsilon$ against $\Delta_t$ should be asymptotic to a straight line with slope $m$
that corresponds to the order of convergence in time of the scheme.
We denote by $\epsilon_2$ and  $\epsilon_{\infty}$ the errors in the
$l_2$-norm and $l_{\infty}$-norm, respectively.
We refer to Table~\ref{tab:1} for the order of convergence in time computed for different fixed mesh widths $h\in \{ 0.1, 0.0.025, 0.00625\}$ and time
steps $\Delta_t \in [T_f/30,T_f/90]$. The solution
computed for $\Delta_t=T_f/100$ is considered as reference solution to
compute the errors. The global errors for the splitting behave like $C(\Delta_t)^2$.
We also observe that the constant $C$ only depends weakly on the spatial
mesh widths $h$.

\begin{table}
\caption{Numerical convergence rates in time for $\theta=\frac{1}{2}$}
\label{tab:1}  
\begin{tabular}{p{2cm}p{3.9cm}p{0.5cm}p{3.9cm}}
\hline\noalign{\smallskip}
Scheme & $l_2$-error convergence rate      && $l_{\infty}$-error convergence rate \\
	& $h=0.1$ \ \ $h=0.025$ \ \ $h=0.00625$ && $h=0.1$ \ \ $h=0.025$ \ \ $h=0.00625$\\
\noalign{\smallskip}\svhline\noalign{\smallskip}
CDS &   2.2002 \hspace*{0.3cm} 2.1975 \hspace*{0.5cm}  2.1969& &2.1973 \hspace*{0.3cm}  2.1958 \hspace*{0.5cm} 2.1956\\
HO5 &    2.1999 \hspace*{0.3cm} 2.1973 \hspace*{0.5cm}  2.1969 &  & 2.1992 \hspace*{0.3cm} 2.1953 \hspace*{0.5cm} 2.1955\\
HOC & 2.2002 \hspace*{0.3cm}    2.1973 \hspace*{0.5cm}   2.1969&  & 2.2007 \hspace*{0.3cm}  2.1953 \hspace*{0.5cm}  2.1955 \\
\noalign{\smallskip}\hline\noalign{\smallskip}
\end{tabular}
\end{table}

In the following, we study the spatial convergence.
The double-logarithmic plots $\epsilon_2$ and $\epsilon_{\infty}$ against
$h$  give the rates of convergence.
Contrary to the time convergence, the order now depends on the parabolic mesh ratio $\mu
=\Delta_t /\Delta_x^2$, so the numerical tests are performed for a set of different constant values of $\mu$.
For simulations, $\mu$ is
fixed at constant values $\mu\in\{0.4, 0.2, 0.1, 0.005 \}$ while
$\Delta_x = \Delta_y= h\to 0$   
($\Delta_t$ is then given by $\Delta_t = \mu h^2$).
The results for the $l_2$-error are given in
Table~\ref{tab:2} and for the $l_{\infty}$-error in Table~\ref{tab:3}.
The solution computed for $h=0.00625$ is used as
reference solution to compute the errors.

\begin{table}
\caption{Numerical convergence rates  in space of $l_2$-error for fixed $\mu$ as $\Delta_x,\,\Delta_t\to 0$ and $\theta=\frac{1}{2}$}
\label{tab:2}  
\begin{tabular}{p{2cm}p{2cm}p{2cm}p{2cm}p{2cm}}
\hline\noalign{\smallskip}
Scheme & $\mu=0.4$         & $\mu=0.2$       & $\mu=0.1$       &  $\mu=0.05$\\
\noalign{\smallskip}\svhline\noalign{\smallskip}
CDS & 1.7828 &  1.7909 &  1.7821 &  1.7845\\
HO5 &      2.2291     &   2.5188  &   2.8153   & 3.0672 \\
HOC & 2.2685   & 2.5191 & 2.8152 & 3.0671 \\
\noalign{\smallskip}\hline\noalign{\smallskip}
\end{tabular}
\end{table}

\begin{table}
\caption{Numerical convergence rates in space of $l_\infty$-error for fixed $\mu$ as $\Delta_x,\,\Delta_t\to 0$ and $\theta=\frac{1}{2}$}
\label{tab:3}  
\begin{tabular}{p{2cm}p{2cm}p{2cm}p{2cm}p{2cm}}
\hline\noalign{\smallskip}
Scheme & $\mu=0.4$         & $\mu=0.2$       & $\mu=0.1$       &  $\mu=0.05$\\
\noalign{\smallskip}\svhline\noalign{\smallskip}
CDS &  1.7170 & 1.7125 &  1.7040 & 1.7038\\
HO5 & 2.2931    &    2.6166 & 2.9182 &3.1584  \\
HOC  & 2.3175 & 2.6176 & 2.9184 & 3.1584  \\
\noalign{\smallskip}\hline\noalign{\smallskip}
\end{tabular}
\end{table}

{\it Remark:}  The choice of the parameter $\theta$ is discussed in \cite{Hund02}. However, for the convergence rates,
$\theta$ seems to have little influence. For example, for the scheme HO5 with  $\theta={1}/{2} +\sqrt{3}/{6}$
we obtain very similar results as shown in Table~\ref{tab:4}.
\begin{table}
\caption{Numerical convergence rates in space for HO5 for fixed $\mu$ as $\Delta_x,\,\Delta_t\to 0$ and $\theta=\frac{1}{2}+ \frac{\sqrt{3}}{6}$}
\label{tab:4}  
\begin{tabular}{p{2cm}p{2cm}p{2cm}p{2cm}p{2cm}}
\hline\noalign{\smallskip}
 & $\mu=0.4$         & $\mu=0.2$       & $\mu=0.1$       &  $\mu=0.05$\\
\noalign{\smallskip}\svhline\noalign{\smallskip}
$l_2$ rate &   2.2310        &2.5186    &    2.8152
 &3.0671 \\
$l_{\infty}$ rate &  2.2938  &   2.6164 &2.9181 &3.1584 \\
\noalign{\smallskip}\hline\noalign{\smallskip}
\end{tabular}
\end{table}

\subsection{Dirichlet Boundary conditions}

In this section we only consider the HOC scheme 
which presents more interesting properties than the other schemes.
Indeed, compared to CDS, its accuracy is larger and compared to HO5, 
no specific treatment at the boundaries is required for the uni-directional terms $F_1$,
$F_2$, the compact scheme is optimal in this respect. 
A particular treatment is necessary when ghost points appear in the
explicit approximation of the mixed term $F_0$.
To preserve the global performance, the accuracy of the approximation near
the boundary conditions has to be sufficiently high. We have used 
a sixth-order approximation in one direction (although lower order
may also be used \cite{GusBC}).
For example, for $u_{0,j}$ on the boundary, at a ghost point $u_{-1,j}$ we impose
$$u_{-1,j}= 5 u_{0,j} -10 u_{1,j}+ 10 u_{2,j} -5 u_{3,j}+u_{4,j}.$$

\begin{figure}
\includegraphics[scale=0.3]{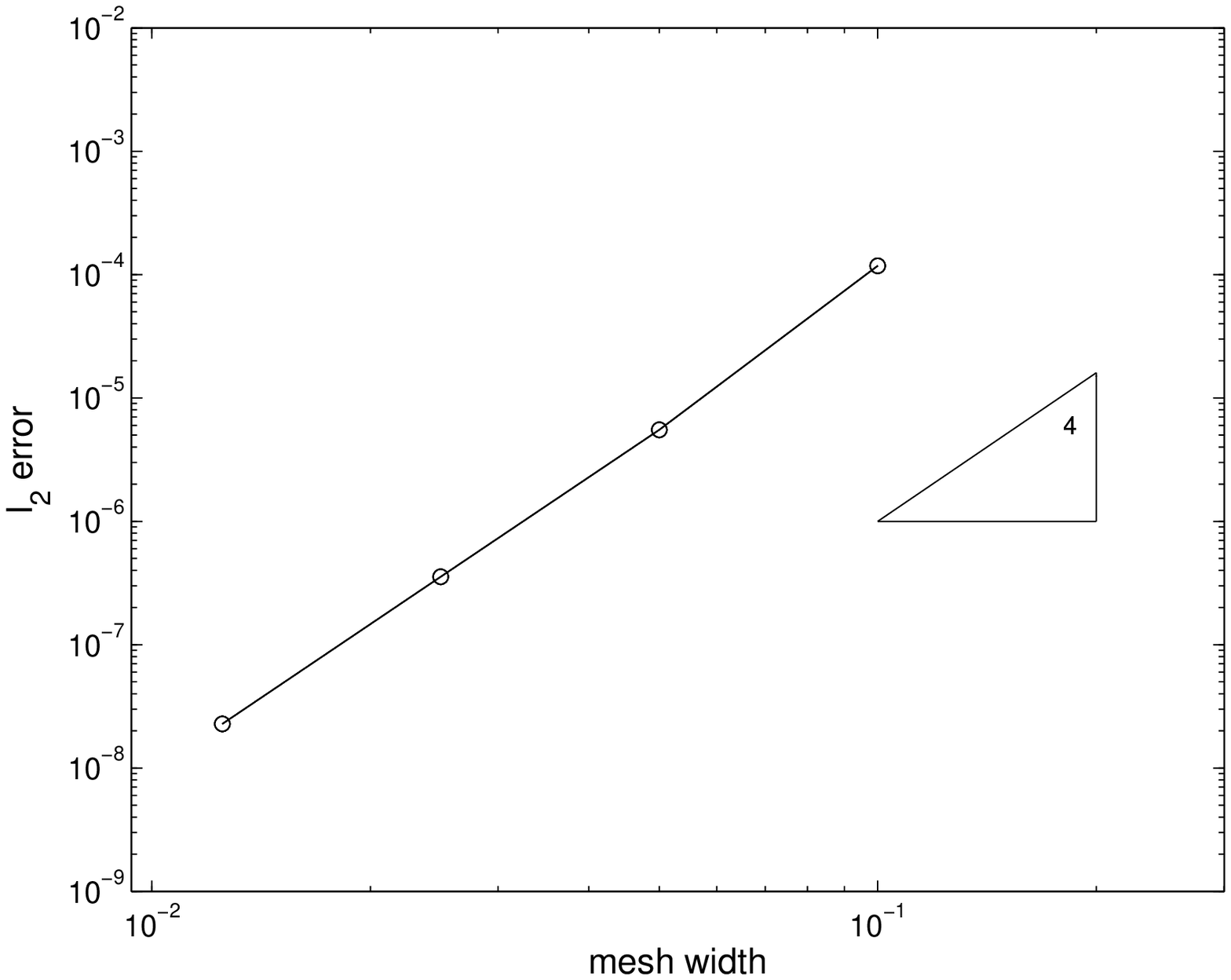} \includegraphics[scale=0.3]{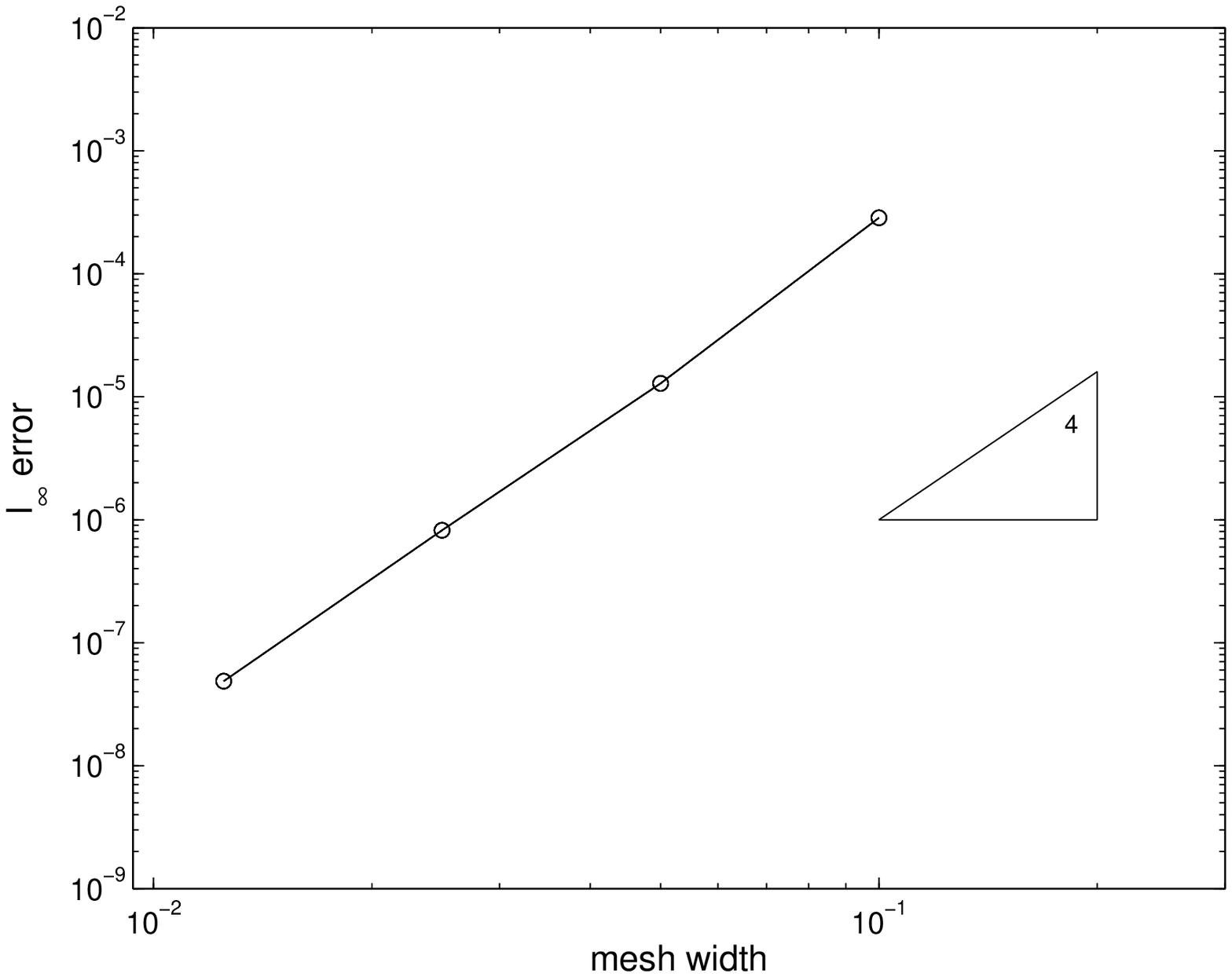}
\caption{Numerical convergence rate in space for HOC ($\theta=\frac{1}{2}$) and $\mu=0.4$.}
\label{fig:1}       
\end{figure}

For the numerical tests, we consider the problem
$$u_t = \mathrm{div} (D\nabla u) + c \cdot \nabla u+ S$$
on the domain $\Omega=[0,1]\times[0,1]$
where
$$
c=-\left ( \begin{array}{c} 
2\\
3\\
\end{array}
 \right ), \quad D= 0.025 \left ( 
\begin{array}{cc} 
1&2\\
2&4\\
\end{array}
\right ),
$$
and the source term $S$ is determined in such a way that the solution is equal to
$u(x,y,t) =- \frac{1}{t+1}\sin(\pi x)\sin(\pi y).$
The Dirichlet boundary condition and initial condition are deduced from the solution.
To incorporate the source term $S$ in the splitting 
(\ref{eq:HVscheme}), $F$ needs to be replaced by $F+S$. More
specifically, $F(U^{n-1})$ is replaced by  $F(U^{n-1}) + S(t^{n-1})$ and $F(Y^2)$ by $F(Y^2) + S(t^{n})$.
We perform the same numerical experiments as in the previous section. The final time is fixed to $T_f = 0.1$ and
the errors are computed with respect to a reference solution computed
on a fine grid in space ($\Delta_x=\Delta_y=0.00625$). Different meshes in space are considered for $\Delta_x=\Delta_y=h$ and
$h\in\{0.1, 0.05, 0.025, 0.0125\}$.
For $\mu=0.4$ the double-logarithmic plots $\epsilon_2$ and $\epsilon_{\infty}$ against
$h$  are given in Figure~\ref{fig:1}.

The results of several numerical tests are reported in Table~\ref{tab:5} for fixed parabolic mesh ratio $\mu =\Delta_t /\Delta_x^2$ while $\Delta_x,\,\Delta_t\to 0$.
In all situations, the new HOC scheme shows a good performance with
fourth-order convergence rates in space, independent of the
parabolic mesh ratio $\mu$.

\begin{table}
\caption{Numerical convergence rates of $l_2$-error and
  $l_{\infty}$-error for HOC ($\theta=\frac{1}{2}$) for different
  constant values of $\mu$ (Dirichlet boundary conditions).}
\label{tab:5}  
\begin{tabular}{p{2cm}p{2cm}p{2cm}p{2cm}p{2cm}}
\hline\noalign{\smallskip}
  & $\mu=0.4$         & $\mu=0.2$       & $\mu=0.1$       &  $\mu=0.05$\\
\noalign{\smallskip}\svhline\noalign{\smallskip}
$l_2$ rate      &  4.0971  &  4.1875 & 4.2129    & 4.2196\\
$l_{\infty}$ rate   & 4.1530 &  4.2372  &  4.2717 & 4.2806\\
\noalign{\smallskip}\hline\noalign{\smallskip}
\end{tabular}
\end{table}


\section{Conclusion}

We have presented new high-order Alternating Direction Implicit (ADI) schemes
for the numerical 
solution of initial-boundary value problems for convection-diffusion equations with mixed derivative terms.
Using the unconditionally stable ADI scheme from \cite{Hund02}
we have proposed different spatial discretizations which lead
to schemes which are fourth-order accurate in space and second-order
accurate in time.

We have performed a numerical convergence analysis with periodic and
Dirichlet boundary conditions where high-order convergence is observed. 
In some cases, the order depends on the
parabolic mesh ratio. More detailed discussions of these schemes including  
this dependence and a stability
analysis will be presented in a forthcoming paper.



\end{document}